\documentclass{amsart}
\usepackage{amsmath, amssymb}
\newtheorem{Theorem}{Theorem}[section]
\newtheorem{Cor}[Theorem]{Corollary}
\newtheorem{Lemma}[Theorem]{Lemma}
\newtheorem{Proposition}[Theorem]{Proposition}
\theoremstyle{definition}
\newtheorem{Definition}[Theorem]{Definition}
\theoremstyle{remark}
\newtheorem{rem}[Theorem]{Remark}
\numberwithin{equation}{section}
\newtheorem{ex}[Theorem]{Example}
\numberwithin{equation}{section}

\newcommand{\R}{\mathbb R}
\newcommand{\N}{\mathbb N}
\newcommand{\Z}{\mathbb Z}

\newcommand{\C}{\mathcal{C}}
\newcommand{\F}{\mathcal{F}}
\newcommand{\p}{\mathcal{P}}

\begin{document}

\title{Remarks on the geometry and the topology of the loop spaces $H^{s}(S^1, N),$ for $s\leq 1/2.$}
\author{Jean-Pierre Magnot}

\address{Lyc\'ee Jeanne d'Arc, Avenue de Grande Bretagne, F-63000 Clermont-Ferrand}%
\email{jean-pierr.magnot@ac-clermont.fr}

\begin{abstract}
We first show that, for a fixed locally compact manifold $N,$ the space $L^2(S^1,N)$ has not the homotopy type of the classical loop space $C^\infty(S^1,N),$ by two theorems:

- the inclusion $C^\infty(S^1,N) \subset L^2(S^1,N)$ is null homotopic if $N $ is connected,

- the space $L^2(S^1,N)$ is contractible if $N$ is compact and connected. 

Then, we show that the spaces $H^s(S^1,N)$ carry a natural structure of Fr\"olicher space, equipped with a 
Riemannian metric, which motivates the definition of Riemannian diffeological space.  
\end{abstract}

\maketitle

\vskip 12pt


\vskip 12pt


\section*{Acknowledgements}
The author is indebted to several participants of the workshop `` diff\'eologies etc..." held in june 2014 in Aix-en-Provence, especially to Paolo Giordano, Seth Wolbert, Jordan Watts, Enxin Wu, for a stimulating discussion on the topic of mapping spaces during lunch, to Augustin Batubenge, for his questions and his great capacity to listen, and to Patrick Iglesias-Zemmour, for his kind invitation to join the meeting and for his stimulating talk. The author would like also to thank very warmly Akira Asada, for exchanges by email around his previous works, which influenced the presentation of part of this paper, and Ekatarina Pervova for email exchanges on possible definitions of Riemannian diffeological spaces.

\tableofcontents

\section*{Introduction}
The loop spaces $H^s(S^1,N),$ for $s>1/2$ are well-known Hilbert manifolds since \cite{Ee}. But very often geometry and analysis stops for $s\leq 1/2,$ because the classical construction of a smooth atlas on these spaces requires an inclusion into the space of continuous loops $C^0(S^1,N), $ via Sobolev embedding theorems. The same holds for loop groups $H^s(S^1,G),$ see e.g. \cite{PS}, where one can read also that, for $s = 1/2,$ most loops $\gamma \in   H^{1/2}(S^1,G) - C^0(S^1,G)$ are not easy to study. One can extend this remark to $s \leq 1/2.$

The aim of this paper is to give a first approach of some topological properties of some of these spaces for $s \leq 1/2,$ and propose an adapted geometric setting. 
In a first part of the paper (section \ref{L2}), we show that there is no homotopy equivalence between $L^2(S^1,N)$ and $H^s(S^1,N)$ for $s>1/2,$
which furnishes a great contrast with the know situation: for $s>k>1/2,$ the inclusion $H^s(S^1,N) \subset H^k(S^1,N)$ is a homotopy equivalence \cite{Pa,Ee,Br}. 

Motivated by the fact that mathematical literature often use weak Sobolev metrics on $C^\infty(S^1,N),$ especially $H^{1/2}$ and $L^2-$metrics (see e.g.\cite{F2,PS,Wu}), 
a natural question is the geometric setting that would enable to discuss with the topologico-geometric properties of the full spaces
$H^s(S^1,N)$ for $s\leq 1/2.$ We then need to find a setting that describes finer structures than the topology, and which enables techniques of differential
geometry. We choose here the setting of Fr\"olicher spaces, which can be seen  as a particular case of diffeological space \cite{Ma2006-3,BIgKWa2014,Wa},
and we develop for the needs of the example of loop spaces the notion of Riemannian diffeological space. 
As a final remark (section \ref{sympl}), we show that the canonical (weak, $H^{1/2}$) symplectic form on the based loop space naturally extends to the (full) based loop $H^{1/2}_0(S^1,N),$
where as the K\"ahler form of the based loop group does not have the same properties.       

 \section{Preliminaries on loop groups and loop spaces}
Let $I = [0;1].$ We note by $(f_n)_{n\in \mathbb{Z}}$ the Fourier coefficients of
any smooth map $f.$ Recall that, for $s \in \mathbb{R}$, 
the space $H^s(I,\mathbb{C})$ 
is the completion of $C^\infty(I,\mathbb{C})$ for 
the norm $||.||_s$ defined by
$$||f||_s^2 = \sum_{n\in\mathbb{Z}}(1+|n|)^{2s}|f_n|^2 = \int_{S^1}  (1 + \Delta^{1/2})^{2s}(f) . \bar{f}  , $$ 
where $\Delta = -\frac{1}{4\pi^2}\frac{d^2}{dx^2}$ is the standard Laplacian, and $S^1 = \R/\Z.$
The same construction holds replacing $\mathbb{C}$ by an algebra $\mathcal{M}$ of matrices with complex coefficients, with the hermitian product of matrices $(A,B) \mapsto tr(AB^*).$
If there is no possible confusion, we note this matrix norm by $||.||$ or by $||.||_{\mathcal{M}}$ if necessary.
Let $N$ be a smooth connected manifold, with Riemannian embedding into $\mathcal{M}.$
We can assume that the $0- $matrix, noted $0,$ is in $N$ with no loss of generality since the space of Riemannian embeddings from $N$ to $\mathcal{M}$ is translation invariant in $\mathcal{M}.$
The loop space $C^\infty(S^1,N)$ is a smooth Fr\'echet manifold (see \cite{Ee,Br} for details). The submanifold of based loops $C^\infty(S^1,N)$ is here identified with loops $\gamma \in C^\infty(S^1,N)$ such that $\gamma(0) = 0.$ 
Let us now consider a compact connected Lie group $G$ of matrices. We note by $C^\infty(S^1,G)$, resp.  $C^\infty_0(S^1,G)$, the group of smooth loops, resp. the group of based smooth loops $\gamma$ such that $\gamma(0) = \gamma(1) = e_G.$  (When dealing about based loop groups, the chosen basepoint is the identity matrix $Id$ for trivial necessities of compatibility with the group multiplication) 

\begin{Definition}
We define $H^s(S^1,G),$ resp. $H^s_0(S^1,G),$ as the adherence of  $C^\infty(S^1,G)$, resp. $C^\infty_0(S^1,G)$, in $H^s(S^1;\mathcal{M}).$ 
\end{Definition}

For $s>1/2,$  it is well-known, that $H^s(S^1,G)$ is a Hilbert Lie group. The key tool is the smooth inclusion $H^s(S^1,\mathcal{M}) \subset C^0(S^1,\mathcal{M}),$
which enables to define charts via tubular neighborhoods, and to define the group multiplication and the group inversion pointwise by the smoothness of the evaluation maps,
see the historical paper \cite{Ee} for details, see also \cite{PS} for an exposition centered on loop groups. 
The biggest Sobolev order where this fails is $s=1/2.$
For $s>1/2,$
\begin{enumerate}
\item the norm $||.||_s$ induces a (strong) scalar product $<.,.>_s$ on $H^s(S^1,\mathfrak{g}),$ which induces a left invariant metric on $TH_s(S^1,G).$
\item if $1/2 < s$, the $H^k-$scalar product $<.,.>_k$ induces a weak Riemannian metric on $TH_s(S^1,G),$ but the $H^k$-geodesic distance is non vanishing on $H^s(S^1,G),$ for $k<s.$
\end{enumerate}

The motivation of this last remark can be found in recent works \cite{BBHM2012,BBHM2013,BBM2013,MichMum} where are given some examples of weak Sobolev $H^s$ metrics  on manifolds of mappings with vanishing geodesic distance. 
\section{The case $s=0:$ On the homotopy type of $L^2(S^1,N)$} \label{L2}

Let us now analyze $L^2(S^1,N)$ when $N$ is connected.

\begin{Lemma} \label{Mol}
Let $ T = \{(l,s) \in I^2 | s \leq l \hbox{ and } l > 0.$
There exists a map $\varphi \in C^\infty(T,[0;1])$ such that 
$$\left\{ \begin{array}{l} \forall l, \varphi(l,0) = 0 \\
\forall l, \varphi(l,l) = 1 \\
\forall l, \frac{\partial\varphi}{\partial s} (l,0) = 1 \\
\forall l, \frac{\partial\varphi}{\partial s} (l,l) = 1 \\
\forall l, \forall k>1, \frac{\partial^k\varphi}{\partial s^k} (l,0) = 0 \\
\forall l, \forall k > 1, \frac{\partial^k\varphi}{\partial s^k} (l,l) = 0 \\
\end{array} \right. .$$
One can choose the following map:
$$ \varphi (l,s) = \int_0^s (\phi_l * m_{1/6l})(t)dt$$
where
\begin{eqnarray*}
\phi_l : & \R & \rightarrow \R \\
& t & \mapsto \left\{ \begin{array}{lr} 1 & \hbox{ if } t < l/3 \hbox{ or } t> 2l/3 \\
\frac{(3 - 2l)}{l} & \hbox{ otherwise } \end{array} \right. 
\end{eqnarray*} 
and 
$m_{1/6l} (t)= 6l m(t/6l),$ with $m$ a standard mollifier with $[-1;1]$ support. 
\end{Lemma}

\noindent
\textbf{Proof.}
The solution given fulfills the conditions required by classical results of analysis. \qed

\vskip 12pt
Notice that with such a function $\varphi, $ we have that $\frac{\partial\varphi}{\partial s} (l,s)\leq 3/l,$ and that $\frac{\partial\varphi}{\partial l} (l,s)\leq M(s)\in \R_+^*,$ for any fixed $0<s\leq 1.$   For a sequence $(s_n)$ such that $s_n \rightarrow 0,$ the sequence $M(s_n)$ is unbounded, were as one can take $M(0)=1.$ These properties are necessary for the proofs of the rest of the section.  

\begin{Theorem} \label{pc}
$L^2_0(S^1,N) = L^2(S^1,N)$ and any loop $\gamma$ in $C^\infty(S^1,N)$ is path-connected  to a null-homotopic piecewise smooth loop in $L^2.$
\end{Theorem}

\noindent
\textbf{Proof.} We now assume that $0 \in N. $ Let $\gamma\in C^\infty(S^1,N)$ such that $\gamma(0) = x \neq 0.$
Let $\tau : [0;1]\rightarrow N$ be a null-homotopic smooth loop such that $\tau(0) = Id,$  $\tau(1/2)=x,$ $\tau(1) = 0,$ $\dot{\tau}(1/2) = \dot\gamma(0)$ in $T_xN.$ Such a loop exists considering a neighborhood in $N$ of a smooth path starting from $0$ and finishing at $x.$
Let $ T = \{(l,s) \in I^2 | s \leq l \hbox{ and } l > 0.$
Let $\varphi \in C^\infty(T,[0;1])$ such that 
$$\left\{ \begin{array}{l} \forall l, \varphi(l,0) = 0 \\
\forall l, \varphi(l,l) = 1 \\
\forall l, \frac{\partial\varphi}{\partial s} (l,0) = 1 \\
\forall l, \frac{\partial\varphi}{\partial s} (l,l) = 1 \\
\forall l, \forall k>1, \frac{\partial^k\varphi}{\partial s^k} (l,0) = 0 \\
\forall l, \forall k > 1, \frac{\partial^k\varphi}{\partial s^k} (l,l) = 0 \\
\end{array} \right. $$  
With the example given in Lemma \ref{Mol}, we have that $$max_{s \in [0;1]} \frac{\partial \varphi}{\partial s} (l,s)\leq 3/l.$$
Let us consider the family of piecewise smooth paths $h \in [0;1] \mapsto \gamma_h$
such that $$ \gamma_h(s) = \left\{ \begin{array}{lcr} \tau \circ \varphi(h/2,s) & \hbox{ if } & s \leq h/2 \\
\gamma  \circ \varphi(1-h, s - h/2 ) & \hbox{ if } & h/2 < s < 1-h/2 \\
 \tau \circ \varphi(h/2,1- s) & \hbox{ if } &  s \geq 1-  h/2 \end{array} \right. .$$
One can check that this is in fact a smooth path, considering the Taylor series at the connecting points $s=h/2,$ $s=1-h/2$ and $s=0.$
Let us take the limit when $h \rightarrow 0.$ 
\begin{eqnarray*}
|| \gamma_h - \gamma ||_{L^2(S^1,\mathcal{M}}^2 & =  & \int_0^{h/2} ||\gamma_h(s) - \gamma(s)||_{\mathcal{M}}^2 ds + \\ &&\int_{h/2}^{1-h/2} ||\gamma_h(s) - \gamma(s)||_{\mathcal{M}}^2 ds + \\ && \int_{1-h/2}^{h/2} ||\gamma_h(s) - \gamma(s)||_{\mathcal{M}}^2ds \\
& =  & \int_0^{h/2} ||\tau \circ \varphi(h/2,s) - \gamma(s)||_{\mathcal{M}}^2 ds + \\ && \int_{h/2}^{1-h/2} ||\gamma \circ \varphi(1-h, s - h/2 ) - \gamma(s)||_{\mathcal{M}}^2 ds +\\
 && \int_{1-h/2}^{h/2} ||\tau \circ \varphi(h/2,1- s) - \gamma(s)||_{\mathcal{M}}^2ds \\
& \leq & \frac{ h  (sup_{s \in [0;1]} ||\tau(s)|| + sup_{s \in [0;1]} ||\gamma(s)||)^2}{2} + \\ && (1-h). \left( hM(1- h)sup_{s \in [0;1]} ||\dot \gamma(s)|| + \frac{3h sup_{s \in [0;1]} ||\dot \gamma(s)||}{2(1-h)}\right)^2 + \\ && \frac{ h  (sup_{s \in [0;1]} ||\tau(s)|| + sup_{s \in [0;1]} ||\gamma(s)||)^2}{2}
\end{eqnarray*}
So that, $$\lim_{h \rightarrow 0} \gamma_h = \gamma,$$
which shows that $\gamma$ is in the $L^2-$ closure of $C^\infty_0(S^1,N).$ Thus we get that $L^2(S^1,N) = L^2_0(S^1,N).$
On the other hand, when we take the $L^2-$limit of $\gamma_h$  when $h \rightarrow 1,$ we get with the same techniques: 
$$ \lim_{h \rightarrow 1} \gamma_h = \tau \vee \tau^{-1}.$$
   \qed

Let us now give another result from the techniques described in the proof of last theorem:

\begin{Theorem} \label{0inj}
The natural injection $C^\infty_0(S^1,N) \rightarrow L^2(S^1,N)$ is homotopic to a constant map.
\end{Theorem}

\noindent
\textbf{Proof.} Let $\gamma\in C^\infty_0(S^1,\mathcal{M}).$
Let $$H(s',\gamma)(l)  = \gamma \circ \varphi(1-s';l) \hbox{ for } 0 \leq l \leq 1-s',$$
with $\varphi$ defined by Lemma \ref{Mol}
that we extend by the constant path on $[1-s;1].$
$H(0,\gamma)= \gamma$ and $H(1,\gamma)$ is the constant loop. 
 For $0<s'<1$ $H(s,\gamma) $ is a piecewise smooth loop, with only one angular point at $l = 1-s',$ and hence is in $L^2.$ We have now to get a continuity property for the map $ s' \mapsto H(s',\gamma).$
Let $s,t$ such that $0  \leq s < t \leq 1$ with $t-s < \frac{s}{6} < \frac{t}{6}.$
We are using in the sequel the following majorations: 

$\bullet$ the classical estimate of the Lipschitz constant of a $C^1$ path: $||\gamma(s) - \gamma(t) || \leq (t-s) max_{s\leq l \leq t}||\dot \gamma (l)||$,

$\bullet$ and since $S^1$ is compact, $ max_{s\in [0;1]}|| \gamma (s)|| = k_\gamma \leq +\infty$

which implies

$\bullet$ on one hand, since $t-s < \frac{s}{6} < \frac{t}{6},$ for $l \in [1-t;1-s],$ $\frac{\partial \varphi}{\partial l} (1-s, l) = 1$, which implies: 
\begin{eqnarray*} \int_{1-t}^{1-s}  ||\gamma \circ \varphi(1-s;l) - \gamma(1)||^2 dl & \leq & max_{0<l<1}  ||\dot \gamma(l)||^2 \int_{1-t}^{1-s}|\varphi(1-s,l) - \varphi(1-s,1-s)|^2 dl\\
& \leq & (t-s)max_{0<l<1}  ||\dot \gamma(l)||^2 \end{eqnarray*}

$\bullet$ and on the other hand, setting $M(l)$ the Lipschitz constant of the map $\phi(.,l),$ : 
$$|\varphi(1-s,l) - \varphi(1-t,l)| \leq M(l)(t-s) \Rightarrow ||\gamma \circ \varphi(1-s,l) - \gamma\circ \varphi(1-t,s)|| \leq M(l)(t-s) max_{0<l<1}  ||\dot \gamma(l)||.$$
On the interval $[(1-t),(1-s)] \subset [0;1[,$ we have that $M(l)$ is bounded by a constant noted $k_{s,t}.$ 

With these two inequalities, we get:
 \begin{eqnarray*}
 || H(s,\gamma) - H(t,\gamma) ||_{L^2(S^1,\mathcal{M})}^2 & =  & \int_0^{1-t} ||\gamma \circ \varphi(1-s;l) - \gamma\circ \varphi(1-t;l)||^2 dl+
\\
&& \int_{1-t}^{1-s}  ||\gamma \circ \varphi(1-s;l) - \gamma(1)||^2 dl \\ && + \int_{1-s}^1 ||\gamma(1) - \gamma(1)||^2 dl \\
& \leq &  k_{s,t}(1-t)(t-s) max_{0<l<1}  ||\dot \gamma(l)||^2 \\ &&+ (t-s) max_{0<l<1}  ||\dot \gamma(l)||^2 + 0
 \end{eqnarray*}

Hence, for a fixed smooth loop $\gamma,$ $s' \mapsto H(s',\gamma) \in C^0(]0;1];L^2(S^1,\mathcal{M})).$
We need to show continuity at $s'=0.$ We get the following inequalities:
 \begin{eqnarray*}
 || H(0,\gamma) - H(t,\gamma) ||_{L^2(S^1,\mathcal{M})}^2 & =  &  \int_{1-t}^{1}  ||\gamma \circ \varphi(1-t;l) - \gamma(1)||^2 dl  \\
& \leq &  t max_{0<l<1}  ||\dot \gamma(l)||^2 
 \end{eqnarray*}
which completes the continuity in the first parameter. We remark that, for fixed $\gamma,$ $H$ is not Lipschitz in the first parameter, since $k_{s,t}$ is not bounded for $(s,t)$ in the neighborhood of $(0;0).$

\vskip 12pt
We now have to show that the map $\gamma \mapsto H(s,\gamma)$ is continuous for the $L^2-$topology. 
For this, we only have to remark the change of coordinates formula:
$$ || \gamma - \tau ||_{L^2(S^1,\mathcal{M})}^2 = \int_0^{1-s}||\gamma \circ \varphi(1-s,l) - \tau\circ \varphi(1-s,l)||^2 \partial_l\varphi(1-s,l) dl.$$
Since $\partial_l\varphi(1-s,l)>1$ for $0 \leq l \leq 1-s,$ we get 
\begin{eqnarray*}
||H(s,\gamma) - H(s,\tau)||_{L^2(S^1,\mathcal{M})}^2 & = & \int_0^1 ||\gamma \circ \varphi(1-s,l) - \tau\circ \varphi(1-s,l)||^2 dl \\
& = & \int_0^{1-s} ||\gamma \circ \varphi(1-s,l) - \tau\circ \varphi(1-s,l)||^2 dl \\
&\leq &  \int_0^{1-s} ||\gamma \circ \varphi(1-s,l) - \tau\circ \varphi(1-s,l)||^2\partial_l\varphi(1-s,l) dl \\
& \leq & || \gamma - \tau ||_{L^2(S^1,\mathcal{M})}^2
\end{eqnarray*}
 Hence the map $\gamma \mapsto H(s,\gamma)$ is $1-$Lipschitz.  \qed

\begin{Cor}
Let $k >1/2.$ The canonical inclusion $i : H^k_0(S^1,N) \rightarrow L^2(S^1,N)$ induces a $0-$map $i_*:H_*(H^k_0(S^1,N)) \rightarrow H_*(L^2(S^1,N))$ and $i_*:\pi_*(H^k_0(S^1,N)) \rightarrow \pi_*(L^2(S^1,N)).$
\end{Cor}

\noindent
\textbf{Proof.}
If $k>1/2,$ the canonical inclusion $C^\infty_0(S^1,N) \rightarrow H^k_0(S^1,N)$ is a well-known homotopy equivalence between smooth manifolds. 
So that, by Theorem \ref{0inj}, we get the result. \qed

\vskip 12pt
We now finish this section with the following result:

\begin{Theorem} \label{contr}
Assume that $N$ is connected and compact. Then the space $L^2(S^1,N)$ is contractible.
\end{Theorem}

\noindent
\textbf{Proof.}
For convenience of the proof, we assume that $0 \in N.$ Let  $H(t,\gamma)(s) = 1_{ s < t } \gamma(s).$ 

$\bullet$ $|| H(t,\gamma) ||_{L^2(S^1,\mathcal{M})} \leq  || 1_{ s < t } ||_{L^\infty(S^1,\mathcal{M})}|| \gamma ||_{L^2(S^1,\mathcal{M})} = || \gamma ||_{L^2(S^1,\mathcal{M})}$ so that $H(t,\gamma)\in L^2(S^1,\mathcal{M}).$ Remarking that $H$ is linear in the second variable, we get that $H(t,.)$ is smooth on $ L^2(S^1,\mathcal{M}).$

$\bullet$ Let $\gamma \in L^2(S^1,N).$ There is a sequence $(\gamma_n)_{n \in \N} \in C^\infty_0(S^1,N)^\N$ that converges to $\gamma.$  

\textit{Claim: The sequence $(H(t,\gamma_n))_{n \in \N}$ is in 
$L^2(S^1,N)$.} 
For this, for fixed $t$ and $n,$ we consider reparametrizations of $\gamma_n$ for $p \in \N^*$ such that $t-1/p < 1:$
$$\delta_p (s) = \left\{\begin{array}{cc} \gamma_n(s) \hbox{ if } & s \leq t \\
 \gamma_n (t + p(s-t)) \hbox{ if } & t < s < t + 1/p \\
0 & \hbox{otherwise} \end{array}\right.$$
We have that the sequence $(\delta_p)$ is in the Sobolev class $H^1, $ and $$||\delta_p -  H(t,\gamma_n) ||_{L^2(S^1,\mathcal{M})} = \left( \int_t^{t+1/p}  \left(\gamma_n (t + p(s-t))\right)^2 ds\right)^{1/2} \leq \frac{||\gamma_n ||_{L^2(S^1,\mathcal{M})}}{p} $$
which shows that $(\delta_p)$ converges to $\gamma_n.$ Since $C^\infty(S^1,N)$ is dense in $H^1(S^1,N)$ \cite{Ee}, we get that   $(H(t,\gamma_n))_{n \in \N}$ is in 
$L^2(S^1,N)$. 

Now, we have  \begin{eqnarray*}||H(t,\gamma_n) - H(t,\gamma) ||_{L^2(S^1,\mathcal{M})} &  \leq &||\gamma_n -\gamma||_{L^2(S^1,\mathcal{M})}
\end{eqnarray*}
So that $H(t,\gamma) \in L^2(S^1,N).$

$\bullet$ Let $\gamma \in L^2(S^1,N).$ For $(t',t') \in [0;1]^2,$ with $t'>t,$ we get \begin{eqnarray*}||H(t',\gamma) - H(t,\gamma) ||_{L^2(S^1,\mathcal{M})} &\leq &||H(t',\gamma) - H(t',\gamma_n) ||_{L^2(S^1,\mathcal{M})} \\ &&
+||H(t',\gamma_n) - H(t,\gamma_n) ||_{L^2(S^1,\mathcal{M})} 
\\ && + ||H(t,\gamma_n) - H(t,\gamma) ||_{L^2(S^1,\mathcal{M})}\\
 &\leq &2||\gamma - \gamma_n ||_{L^2(S^1,\mathcal{M})} \\ &&
+||H(t',\gamma_n) - H(t,\gamma_n) ||_{L^2(S^1,\mathcal{M})} 
\end{eqnarray*}
Now, N is compact so that, $k = sup_{x \in N} ||x||_\mathcal{M} < +\infty.$ For this, we get \begin{eqnarray*}   ||H(t',\gamma_n) - H(t,\gamma_n) ||_{L^2(S^1,\mathcal{M})} & \leq & k ||1_{ [t,t']}||_{L^2(S^1,\mathcal{M})} \\& = & k(t'-t) \end{eqnarray*}
Inserting this last inequality in the previous one, for $\epsilon >0,$ choose $ (t'-t)<\epsilon/3k$ and $n$ such that  $ ||\gamma - \gamma_n ||_{L^2(S^1,\mathcal{M})} < \epsilon/3,$
we get that $||H(t',\gamma) - H(t,\gamma) ||_{L^2(S^1,\mathcal{M})} \leq \epsilon,$ which shows $H$ is continuous in the first variable, and ends the proof. \qed

\noindent
\begin{rem}
The same procedure can be adapted replacing $L^2(S^1,N)$ by $L^2(M,N),$ with $M$ compact manifold. With the same arguments, one can build the homotopy map with a smooth Morse function, and mimick line by line the last proof. This proof will be developped elsewhere, for the sake of unity of the exposition.
\end{rem}

\section{Riemannian metrics and Hausdorff measures on diffeological spaces}

Diffeological spaces and Fr\"olicher spaes will furnish a setting to deal with the differential geometry of the loop spaces $H^s(S^1, N).$
For preliminaries on diffeological spaces and Fr\"olicher spaces, we refer to \cite{Igdiff} and to \cite{FK,KM}. For convenience, the necessary material and a complementary bibliography 
is given in the appendix.
We now describe an extension of some basic structures of Riemannian manifolds to diffeological spaces. 
 \subsection{Riemannian diffeological spaces} 

\begin{Definition} 
Let $(X,\p)$ be a diffeological space.
A \textbf{Riemannian metric} $g$ on $X$ (noted $g \in Met(X)$) is a map 
$$  g: \{p:O_p\rightarrow X\} \in \p \mapsto g_p$$
such that 
\begin{enumerate}
\item $x \in O_p \mapsto g_p(x)$ is a smooth section of the bundle of symmetric bilinear forms on $TO_p$  

\item let $y\in O_p$ and $y'\in O_{p'}$ such that $p(y)=p'(y').$ If $(X_1,X_2)$ is a couple of germs of paths in 
$Im(p)\cap Im(p'),$ if there exists two systems of $2-$vectors $(Y_1,Y_2)\in (T_yO_p)^2$ and $(Y'_1,Y'_2)\in (T_{y'}O_{p'})^2,$ if $p_*(Y_1,Y_2)=p'_*(Y'_1,Y'_2)=(X_1,X_2),$
$$ g_p(Y_1,Y_2) = g_{p'}(Y'_1,Y'_2).$$

\item \label{c3} for each non zero germ of smooth path $Y,$ $$g(Y,Y)>0.$$ \end{enumerate}
\noindent
$(X,\p,g)$ is a \textbf{Riemannian diffeological space} if $g$ is a metric on $(X,\p).$ If condition \ref{c3} is not everywhere fulfilled, we call it \textbf{pseudo-Riemannian diffeological space.}
\end{Definition}
For any germ of path $X$ we note $||X|| = \sqrt{g(X,X)}.$
\begin{Definition}
We call \textbf{arc length} the map  $L : C^\infty([0;1], X) \rightarrow \mathbb{R}_+$ defined by $$L(\gamma) = \int_0^1 ||\dot\gamma(t)|| dt.$$
Let $(x, y) \in X^2.$ We define $$ d_g(x,y) = \inf \left\{ L(\gamma) | \gamma(0) = x \wedge \gamma(1)=y\right\}$$
and we call \textbf{Riemannian pseudo-distance} the map $d : X \times X \rightarrow \R_+$ that we have just described. 
\end{Definition}

The following proposition justifies the terms ``pseudo-distance'': 
\begin{Proposition}
\begin{enumerate}
\item \label{d1} $\forall x \in X, d_g(x,x) = 0.$
\item \label{d2} $\forall (x,y) \in X^2, d_g(x,y) = d_g(y,x)$
\item \label{d3} $\forall (x,y,z) \in X^3, d_g(x,z) \leq d_g(x,y) + d_g(y,z).$  
\end{enumerate}
\end{Proposition}

\noindent
\textbf{Proof.} The proofs are standard, let us recall the main arguments. For \ref{d1}, the constant path gives the minimum. For \ref{d2}, the reverse parametrization $t \mapsto 1-t$  defines a transformation 
from the paths from $x$ to $y$ to the paths from $y$ to $x$ under which $L$ is invariant. 
For \ref{d3}, the paths passing by $y$ are only a part of the set of paths from $x$ to $z.$ \qed

One could wonder whether $d$ is a distance or not, i.e. if we have the stronger property: 
$$ \forall (x,y) \in X^2, d_g(x,y)=0 \Leftrightarrow x=y.$$
Unfortunately, it seems to appear in examples arising from infinite dimensional geometry that there can 
have a distance which equals to 0 for $x \neq y.$ This is what is described on e.g. a weak Riemannian metric 
of a space of proper immersions in the work of Michor and Mumford \cite{MichMum}. 
Moreover, the D-topology is not the topology defined by the pseudo-metric $d.$ All these facts, which show 
that the situation on Riemannian diffeological spaces is different from the one known on finite dimensional manifolds, are checked in the following (toy) example.

\begin{ex}
Let $Y = \coprod_{i \in \N^*} \R_i$ where $\R_i$ is the $i-$th copy of $\R,$ equipped with its natural scalar product.
Let $\mathcal{R}$ be the equivalence relation $$x_i \mathcal{R} x_j \Leftrightarrow  \left\{ \begin{array}{l}(x_i \notin ]0;\frac{1}{i}[ \wedge x_j \notin ]0;\frac{1}{j}[) \Rightarrow \left\{ \begin{array}{lr} x_i = x_j & \hbox{if } x_i\leq 0 \\
x_i + 1- \frac{1}{i}  = y_j + 1 - \frac{1}{j} & \hbox{ if } x_i \geq \frac{1}{i} \end{array} \right. \\ (x_i \in ]0;\frac{1}{i}[ \vee x_j \in ]0;\frac{1}{j}[) \Rightarrow i = j \wedge x_i = x_j  \end{array} \right. $$
Let $X = Y / \mathcal{R}.$ This is a 1-dimensional Riemannian diffeological space. Let $\bar{0}$ be the class of $0 \in\R_1,$ and let $\bar{1}$ 
be the class of $1 \in \R_1.$  Then $d_g(\bar{0}, \bar{1}) = 0.$ This shows that $d_g$ is not a distance on $X.$  In the $D-$topology, $\bar{0}$ and $\bar{1}$ respectively have the following disconnected neighborhoods: 
$$U_{\bar{0}} = \left\{ \bar{x_i} | x_i < \frac{1}{2i} \right\}$$
and $$U_{\bar{1}} =\left\{\bar{x_i} | x_i > \frac{1}{2i} \right\}.$$
This shows that $d$ does not define the $D-$topology.  
\end{ex}

This leads to the following definition:

\begin{Definition}
A \textbf{Fr\"olicher-strong} Riemannian metric is a Riemannian metric $g$ for which $d$ is a distance.
\end{Definition}

We have to notice that this notion is not exactly the same as the notion of strong Riemannian metric on a Hilbert vector bundle.

\subsection{Volume and diffeologies}

On a (finite dimensional) Riemannian manifold $M,$ the notion of Riemannian volume is related to the dimension of the manifold, and to the notion of volume form $\omega_M.$ 
On one hand, if the Riemannian manifold $M$ is viewed now as a Fr\"olicher space, with underlying diffeology $\p$, we have that 
$$\forall f \in \p, f^*\omega_M = 0 \Leftrightarrow \hbox{ Dimension of } f < dim M.$$
On the other hand, if $f \in \p$ is an embedding  $O \rightarrow M,$ and if $O$ is an open domain of dimension $p$, $f(O)$ is a submanifold of $M$ and it can be equipped with the $p-$dimensional Hausdorff 
measure induced by the geodesic distance on $O,$ and we have:

\begin{Proposition}
Assume that $O \subset \R^m$ is a $n-$dimensional submanifold. The Hausdorff dimension of $O$ is $n$ and, for any relatively open subset $U \subset O,$
if $\mathcal{H}^n$ is the $n-$dimensional Hausdorff measure in $\R^m,$ $$\mathcal{H}^n(U) = \int_U \omega_O .$$
\end{Proposition}

We remark that, given a chart $O$ on $M,$ $O$ is equipped with the standard Lebesgue volume $d\lambda = \omega_O,$ and $M$ is equipped with the Riemannian volume $\omega_M,$
on $O,$ $$ \omega_M = \sqrt{det g} . d\lambda $$
 and hence, if $U$ is a n-dimensional submanifold of $O,$ noting $i_u$ the canonical injection $U \rightarrow O,$ we can define $$ \mathcal{H}^n_M(U) = \int_U \sqrt{det(i^*g)} d\mathcal{H}^n .$$
This corresponds to the $n$ dimensional Hausdorff measure with respect to $M$ viewed as a metric space. As a consequence, a Riemannian manifold does not only carry one volume measure, but a family of measures on the plots of its diffeology. If $f$ is a $n$dimensional plot of the diffeology of $M,$ we define $\mathcal{H}_f = f^* \mathcal{H}^n_M.$
This property is stable under composition of plots, reparametrization, gluing. This leads to the following on Riemannian diffeological spaces:

\begin{Definition}
Let $(X,\p)$ be a diffeological space.
The \textbf{Hausdorff diffeological volume} associated to a Riemannian metric is the collection $\left\{\mathcal{H}_{p}; p \in \p \right\}$ of $dim(D(p))$ Hausdorff measures on the domains $D(p).$
\end{Definition}

Let us remark that:
\begin{itemize}
\item if $det(p^*g) > 0$ on $D(p),$ the domain of $p$, then $\mathcal{H}_p$ is $dim(D(p))$-dimensional Hausdorff measures on $D(p)$ induced by the Riemannian distance on $D(p) \subset X.$ 
\item for any $(p,p') \in \p^2$ of same dimension, if  $p'= f \circ p,$ $\mathcal{H}_p' = f^* \mathcal{H}_p.$
\item If there exists $x \in D(p)$ such that $det(p^*g_x)=0,$ the definition of the Hausdorff metric via the (pseudo)-distance on $D(p)$ remains valid \cite{Fed}. 
 \end{itemize}

We have here to remark that the Riemannian metric needs not to be Fr\"olicher-strong, because this is the induced Riemannian metric on each $D(p)$ which defines the Hausdorff measure, one should say on the (classical) Riemannian manifold $D(p).$ 
\subsection{On $\infty-p$ forms and volume forms}
This section is based on ideas from A. Asada \cite{As2004,As2005} adapted to the context of a diffeological space . 
\begin{Definition}
Let $(X,\p)$ be a Riemannian diffeological space. An \textbf{orientable plot} on $X$ is a plot $p\in \p$, $dim(p)=n,$ such that there exists a n-form $\omega_n\in \Omega^n(X,\R)$ such that $p_*h_p =\omega_n,$ where $h_p$ is a $n-$form on $D(p),$ that induces the Hausdorff measure $\mathcal{H}_p.$ The space of orientable plots of the diffeology $\p$ is noted $\mathcal{O}(\p).$  
\end{Definition}

\begin{Proposition}
Let $X$ be a smooth n-dimensional compact manifold, equipped with its n\'ebuleuse diffeology $\p.$ X is orientable if and only if there exists a n-plot $p \in \p,$ surjective, such that $p \in \mathcal{O}(\p).$
\end{Proposition}

The proof is straightforward, setting a Riemannian metric $g$ on $X,$ and using the exponential map to define the plot $p.$ 

\begin{ex}
Let $X=S^n,$ $n\geq 1.$ Let $P,P'$ be two antipodal points. The mapping $exp_P: T_PS^n \rightarrow S^n$ has an injectivity radius $r = \pi.$ The cut locus is $P'.$
Thus, considering the open ball $B(0,3\pi/2) \subset T_PS^n$ and the plot  $p : exp_P|_{B(0,3\pi/2)},$ we get the construction, even if $p^*\omega_{S^n} = 0$ on $p^{-1}(P')$ (here, $\omega_{S^n}$ is the canonical volume form on $S^n$). 
\end{ex}

\begin{ex} Let $X' = ]0;1[ \times [0;1],$ let $\sim$ be the relation of equivalence on $X'$ defined by $(t,0) \sim (1-t,1)$ and let $X = X'/\sim$ be the (open) M\"obius band. The mapping $ p:]0;1[ \times ]-1/2;3/2[ \rightarrow X'$ defined by $$ p(x,y) = \left\{ \begin{array}{lcr} (x,y+1) & \hbox{if} & x<0 \\ 
(x,y) & \hbox{if} & x \in [0;1] \\
 (x,y-1) & \hbox{if} & x>1 \end{array} \right. $$
is a 2-dimensional plot such that the trace of the canonical Lebesgue measure co\"incides with $\mathcal{H}_p,$ but for which there exists no 2-form $\omega_2$ on $X'$ such that  $p_*\lambda = \omega_2$ where $\lambda$ is the canonical Lebesgue measure on $]0;1[ \times ]-1/2;3/2[$, because $\omega_2$ should be non zero everywhere.
\end{ex}

After these examples, let us turn to the key definition :

\begin{Definition}
Let $(X,\p)$ be a Riemannian diffeological space with set of oriented plots $\mathcal{O}\p.$
We call \textbf{volume form } of $X$ a collection of forms 
$$p \in \mathcal{O}\p \mapsto \omega_p \in \Omega^{dim D(P)}(D(p), \R)$$
such that 

- the form $\omega_p$ defines the $dim(D(p))-$ Hausdorff measure on $D(p)$

- if $p$ and $ p'$ are oriented $n-$plots such that $p' = p \circ f,$ then $\omega_p' = f^* \omega_p $
\end{Definition}

\begin{Definition}
Let $(X,\p)$ be a Riemannian diffeological space with volume form $\omega.$ 
Let $q \in \N.$ 
A $(\infty-q)-$form is a collection
 $$p \in \mathcal{O}\p \mapsto \omega_p \in \Omega^{dim D(P)-q}(D(p), \R)$$
such that there exists a $q-$form $\beta \in \Omega^q(X,\R)$ such that 
$$ \forall p \in \mathcal{O}\p,  \alpha_p \wedge p^*\beta = \omega_p.$$
For the consistency of the definition, if $q > dim(p)$ or if $\omega_p = 0,$ we set $\alpha_p = 0.$

\end{Definition}

With such a definition, a volume form is a $(\infty - 0)-$form.
We note by $\Omega^{(\infty-q)}(X,\R)$ the space of $(\infty-q)-$forms.

\section{$H^{s}(S^1,N)$ as a Riemannian diffeological space}

\subsection{Settings}

\begin{Proposition}
Let $s\leq 1/2.$. Then $H^{s}(S^1,N)$ and $H^{s}_0(S^1,N)$ are Riemannian Fr\"olicher space. The same holds for $N=G.$
\end{Proposition}

\noindent
\textbf{Proof.}
$H^{s}(S^1,\mathcal{M})$ is equipped with its natural underlying structure of Hilbert space, which carries a natural structure of Fr\"olicher space. 
As subsets,  $H^{s}(S^1,N)$ and $H^{s}_0(S^1,N)$ are equipped with the reflexive completion of their trace diffeology. 
So that, they are Fr\"olicher spaces. The natural Hilbert structure on $H^s(S^1,\mathcal{M})$ induces a Riemannian metric on $H^s(S^1,N).$\qed

\vskip 12pt

\begin{Proposition}
$H^s(S^1,N)$ is Fr\"olicher-strong for $s \in \R.$
\end{Proposition}

\noindent
\textbf{Proof.} Let $\gamma$ be a smooth path in $H^s(S^1,N) \subset H^s(S^1,\mathcal{M}).$ Then the length of $\gamma$ is bounded below by $||\gamma(1)-\gamma(0)||_{H^s(S^1, \mathcal{M})}$ \qed 

Let us now give a result for the extension of the multiplication of loop groups. For this, we define the space
$$H^{1/2,+}(S^1,G) = \bigcup_{s>1/2} H^s(S^1,G)$$

\begin{Lemma}
The space $H^{1/2,+}(S^1,G)$ is a Lie group modeled on a locally convex vector space.
\end{Lemma}

\noindent
\textbf{Proof.}
For each $s>1/2,$ we have $H^s(S^1,G) \subset C^0(S^1,G)$ and the usual (exponential) atlas on $H^s(S^1,G)$ is induced by the atlas on $C^0(S^1,G),$ see \cite{Ee,PS} for the details.
Then, following \cite{Glo}, we get the result. \qed

\begin{Proposition} Let $k > 1/2$ and let $s \leq k.$ The natural action $ C^\infty(S^1,G) \times C^\infty(S^1,G) \rightarrow C^\infty(S^1,G)$ 
extends to a smooth action $$ H^k(S^1,G) \times H^s(S^1,G) \rightarrow H^s(S^1,G),$$ and for $s\leq 1/2,$ to a smooth action
$$ H^{1/2,+}(S^1,G) \times H^s(S^1,G) \rightarrow H^s(S^1,G).$$
\end{Proposition}\noindent
\textbf{Proof.} Since $G \subset \mathcal{M},$ with smooth inclusion and trace diffeology, it is sufficient to remark that this theorem is an 
application of the standard ``multiplication theorem' of Sobolev classes, which states that the multiplication, with the coefficients defined as above, 
is a bilinear continuous map. \qed.

\subsection{$H^{1/2}_0(S^1,N)$ and symplectic forms} \label{sympl}

Let $\gamma \in H^{1/2}_0(S^1,\mathcal{M}).$ Then $\dot{\gamma} \in H^{-1/2}(S^1,\mathcal{M})$ The canonical 1-form
$$\theta(X) = \int_{S^1} (\dot{\gamma}, X)$$
defined first for $\gamma \in C^\infty_0(S^1,\mathcal{M}$ and  $X \in C^\infty(S^1,\mathcal{M})$ extends to a 1-forms on $H^{1/2}_0(S^1,\mathcal{M})$ by the pairing of dual spaces $$H^{-1/2} \times H^{1/2}_0 \rightarrow \R.$$
So that, 

\begin{Theorem}
The symplectic 2-form on the based loop space $C^\infty_0(S^1,N)$ defined by 
$$ \omega_N(X,Y) = d\theta(X,Y) = \int_{S^1} \left( \frac{\nabla^N}{ds} X(s), Y(s) \right) ds$$
extends to a 2-form on $H^{1/2}_0(S^1, N).$
\end{Theorem}

This is also the case when $N=G.$ On $C^\infty_0(S^1,G),$ the vector field $\dot{\gamma}$ is not left-invariant. 
We know that, on the based loop group, there is another symplectic form, which is not exact, defined for left-invariant vector fields $X $ and $Y$ by 
$$ \omega_G(X,Y) = \int_{S^1} \left( \frac{d X(s) }{ds}, Y(s) \right) ds.$$

But the 2-form $\omega_G$ does not seem to extend to $H^{1/2}_0(S^1, G)$ because the full space $H^{1/2}_0(S^1, G)$ is not a group. The biggest known group in $H^{1/2}_0(S^1, G)$ is $H^{1/2}_0(S^1, G)\cap C^0(S^1, G),$ see \cite{PS}.

\section{Appendix: Diffeological and Fr\"olicher spaces}

\label{1.} The objects of the category of -finite or infinite- dimensional
smooth manifolds is made of topological spaces $M$ equipped with
a collection of charts called maximal atlas that enables one to make
differentiable calculus. But in examples of projective limits of manifolds,
a differential calculus is needed where as no atlas can be defined.
To circumvent this problem which occurs in various frameworks, several
authors have independently developed some ways to define differentiation
without defining charts. We use here two of them. The first one
is due to Souriau (\cite{Sou}, see e.g. \cite{Igdiff} for a textbook), the second one is due to Fr\"olicher (see \cite{FK}, and e.g. \cite{KM} for an introduction). In this section, we review some basics
on these notions. A (non exhaustive) complementary list of reference is \cite{BIgKWa2014,BN2005,BT2014,CSW2014,CW2014,CN,DN2007-1,DN2007-2,Les,Ma2006-3,Ma2013,Wa}.

\subsection{Diffeological spaces and
Fr\"olicher spaces}
\label{1.1}

\begin{Definition} Let $X$ be a set.

\noindent $\bullet$ A \textbf{parametrization} of dimension $p$ (or $p$-plot)
on $X$ is a map from an open subset $O$ of $\R^{p}$ to $X$.

\noindent $\bullet$ A \textbf{diffeology} on $X$ is a set $\p$
of parametrizations on $X$ such that, for all $p\in\N$,

- any constant map $\R^{p}\rightarrow X$ is in $\p$;

- Let $I$ be an arbitrary set; let $\{f_{i}:O_{i}\rightarrow X\}_{i\in I}$
be a family of maps that extend to a map $f:\bigcup_{i\in I}O_{i}\rightarrow X$.
If $\{f_{i}:O_{i}\rightarrow X\}_{i\in I}\subset\p$, then $f\in\p$.

- (chain rule) Let $f\in\p$, defined on $O\subset\R^{p}$. Let $q\in\N$,
$O'$ an open subset of $\R^{q}$ and $g$ a smooth map (in the usual
sense) from $O'$ to $O$. Then, $f\circ g\in\p$.

\noindent
\vskip 6pt $\bullet$ the parametrizations $p \in \p$ are called the {\bf plots} of the diffeology $\p.$

\noindent\vskip 6pt $\bullet$ If $\p$ is a diffeology $X$, $(X,\p)$ is
called \textbf{diffeological space}.

\noindent Let $(X,\p)$ et $(X',\p')$ be two diffeological spaces,
a map $f:X\rightarrow X'$ is \textbf{differentiable} (=smooth) if
and only if $f\circ\p\subset\p'$. \end{Definition}

\begin{rem} Notice that any diffeological space $(X,\p)$
can be endowed with a natural topology such that all the maps that
are in $\p$ are continuous.This topology  is called the $D-$topology \cite{CSW2014}.
\end{rem}
\begin{rem}
Let $f\in\p$, defined on $O\subset\R^{p}$. we call $p$ the dimension of the plot $f.$ 
\end{rem}

We now introduce Fr\"olicher spaces.

\begin{Definition} $\bullet$ A \textbf{Fr\"olicher} space is a triple
$(X,\F,\C)$ such that

- $\C$ is a set of paths $\R\rightarrow X$,

- A function $f:X\rightarrow\R$ is in $\F$ if and only if for any
$c\in\C$, $f\circ c\in C^{\infty}(\R,\R)$;

- A path $c:\R\rightarrow X$ is in $\C$ (i.e. is a \textbf{contour})
if and only if for any $f\in\F$, $f\circ c\in C^{\infty}(\R,\R)$.

\vskip 5pt $\bullet$ Let $(X,\F,\C)$ et $(X',\F',\C')$ be two
Fr\"olicher spaces, a map $f:X\rightarrow X'$ is \textbf{differentiable}
(=smooth) if and only if $\F'\circ f\circ\C\in C^{\infty}(\R,\R)$.
\end{Definition}

Any family of maps $\F_{g}$ from $X$ to $\R$ generate a Fr\"olicher
structure $(X,\F,\C)$, setting \cite{KM}:

- $\C=\{c:\R\rightarrow X\hbox{ such that }\F_{g}\circ c\subset C^{\infty}(\R,\R)\}$

- $\F=\{f:X\rightarrow\R\hbox{ such that }f\circ\C\subset C^{\infty}(\R,\R)\}.$

A Fr\"olicher space carries a natural topology,
which is the pull-back topology of $\R$ via $\F$, see e.g. \cite{BT2014}. In the case of
a finite dimensional differentiable manifold, the underlying topology
of the Fr\"olicher structure is the same as the manifold topology. In
the infinite dimensional case, these two topologies differ very often.

\vskip 12pt

In the previous settings, we call $X$ a \textbf{differentiable
space}, omitting the structure considered. Notice that the sets of differentiable maps between two differentiable
spaces satisfy the chain rule. Let us now compare these settings:
Let $(X,\F, \mathcal{C})$ be a Fr\"olicher space. We define
$$
\p(\F)=
\coprod_{p\in\N}\{\, f\hbox{ p-
paramatrization on } X; \, \F \circ f \in C^\infty(O,\R) \quad \hbox{(in
the usual sense)}\}.$$
With this construction, we also get a natural diffeology when
$X$ is a Fr\"olicher space, extension of the ``n\'ebuleuse'' diffeology of a manifold \cite{Sou, Igdiff}. In this case, one can easily show the following:
\begin{Proposition} \label{Frodiff} \cite{Ma2006-3} 
Let $(X,\F,\C)$
and $(X',\F',\C')$ be two Fr\"olicher spaces. A map $f:X\rightarrow X'$
is smooth in the sense of Fr\"olicher if and only if it is smooth for
the underlying diffeologies. \end{Proposition}

Thus, we can state in an intuitive but comprehensive way:

\vskip 6pt
\begin{tabular}{ccccc}
smooth manifold  & $\Rightarrow$  & Fr\"olicher space  & $\Rightarrow$  & Diffeological space\tabularnewline
\end{tabular}

\vskip 6pt
The classical properties of a tangent space are not automatically checked, neither in the case of diffeological spaces following \cite{CW2014}, nor in the case of Fr\"olicher spaces following \cite{DN2007-1}.  In these two settings, there are two possible tangent spaces:

-the \textbf{internal} tangent space, defined by the derivatives of smooth paths,

- the \textbf{external} tangent space, made of derivations on $\R-$valued smooth maps. 

The internal tangent space is not necessarily a vector space, were as the external tangent space is. For example, consider the diffeological subspace of $\R^2$ made of the two lines $y=x$ and $y = -x.$ They cross at the origin, the external tangent space is of dimension 2, the internal tangent space is made of two directions that cannot be combined by addition. This is why some authors sometimes call the internal tangent space by ``tangent cone''. A refinement of the internal tangent space has been recently studied in \cite{CW2014}.
\subsection{Fr\"olicher completion of a diffeological space} \label{complet}

We now finish the comparison of the notions of diffeological and Fr\"olicher 
space following mostly \cite{Wa}:

\begin{Theorem}
Let $(X,\p)$ be a diffeological space. There exists a unique Fr\"olicher structure
$(X, \F_\p, \C_\p)$ on $X$ such that for any Fr\"olicher structure $(X,\F,\C)$ on $X,$ these two equivalent conditions are fulfilled:

(i)  the canonical inclusion is smooth in the sense of Fr\"olicher $(X, \F_\p, \C_\p) \rightarrow (X, \F, \C)$

(ii) the canonical inclusion is smooth in the sense of diffeologies $(X,\p) \rightarrow (X, \p(\F)).$ 

\noindent Moreover, $\F_\p$ is generated by the family 
$$\F_0=\lbrace f : X \rightarrow \R \hbox{ smooth for the 
usual diffeology of } \R \rbrace.$$
\end{Theorem}

\noindent
\textbf{Proof.}
Let $(X,\F,\C)$ be a Fr\"olicher structure satisfying \textit{(ii)}. 
Let $p\in P$ of domain $O$. $\F \circ p \in C^\infty(O,\R)$ in the usual sense.  
Hence, if $(X,\F_\p, \C_\p)$ is the Fr\"olicher structure on $X$ generated by the 
set of smooth maps $(X,\p)\rightarrow \R,$ we have two smooth inclusions 
$$ (X,\p) \rightarrow (X,\p(\F_\p)) \hbox{ in the sense of diffeologies }$$
and
$$ (X, \F_\p, \C_\p) \rightarrow (X,\F,\C) \hbox{ in the sense of Fr\"olicher. }$$
Proposition \ref{Frodiff} ends the proof. \quad \qed

\begin{Definition} \cite{Wa}
A \textbf{reflexive} diffeological space is a diffeological space $(X,\p)$ such that $\p = \p(\F_\p).$
\end{Definition}

\begin{Theorem} \cite{BIgKWa2014,Wa}
The category of Fr\"olicher spaces is exactly the category of reflexive diffeological spaces.
\end{Theorem}

This last theorem allows us to make no difference between Fr\"olicher spaces and reflexive diffeological spaces. 
We shall call them Fr\"olicher spaces, even when working with their underlying diffeologies.

\subsection{Push-forward, quotient and trace}

We give here only the results that will be used in the sequel.

\begin{Proposition} \cite{Ma2006-3} Let $(X,\p)$ be a diffeological space,
and let $X'$ be a set. Let $f:X\rightarrow X'$ be a surjective map.
Then, the set \[
f(\p)=\{u\hbox{ such that }u\hbox{ restricts to some maps of the type }f\circ p;p\in\p\}\]
 is a diffeology on $X'$, called the \textbf{push-forward diffeology}
on $X'$ by $f$. \end{Proposition}

Given a subset $X_{0}\subset X$, where $X$ is a Fr\"olicher space
or a diffeological space, we can define on trace structure on $X_{0}$,
induced by $X$.

$\bullet$ If $X$ is equipped with a diffeology $\p$, we can define
a diffeology $\p_{0}$ on $X_{0}$ setting \[
\p_{0}=\lbrace p\in\p\hbox{such that the image of }p\hbox{ is a subset of }X_{0}\rbrace.\]

$\bullet$ If $(X,\F,\C)$ is a Fr\"olicher space, we take as a generating
set of maps $\F_{g}$ on $X_{0}$ the restrictions of the maps $f\in\F$.
In that case, the contours (resp. the induced diffeology) on $X_{0}$
are the contours (resp. the plots) on $X$ which image is a subset
of $X_{0}$.

\subsection{Cartesian products and projective limits}

\begin{Proposition} \label{prod1} Let $(X,\p)$ and $(X',\p')$
be two diffeological spaces. We call \textbf{product diffeology} on
$X\times X'$ the diffeology $\p\times\p'$ made of plots $g:O\rightarrow X\times X'$
that decompose as $g=f\times f'$, where $f:O\rightarrow X\in\p$
and $f':O\rightarrow X'\in\p'$. \end{Proposition}

In the case of a Fr\"olicher space, we derive very easily, compare
with e.g. \cite{KM}: \begin{Proposition} \label{prod2} Let $(X,\F,\C)$
and $(X',\F',\C')$ be two Fr\"olicher spaces, with natural diffeologies
$\p$ and $\p'$ . There is a natural structure of Fr\"olicher space
on $X\times X'$ which contours $\C\times\C'$ are the 1-plots of
$\p\times\p'$. \end{Proposition}

We can even state the same results in the case of infinite products,
in a very trivial way by taking the cartesian products 
of the plots or of the contours.  Let us now give the description of what happens
for projective limits of Fr\"olicher and diffeological spaces.

\begin{Proposition} \label{froproj} Let $\Lambda$ be an infinite
set of indexes, that can be uncoutable.

$\bullet$ Let $\lbrace(X_{\alpha},\p_{\alpha})\rbrace_{\alpha\in\Lambda}$
be a family of diffeological spaces indexed by $\Lambda$ totally
ordered for inclusion,  with $(i_{\beta,\alpha})_{(\alpha,\beta)\in\Lambda^{2}}$
a family of diffeological maps . If $X=\bigcap_{\alpha\in\Lambda}X_{\alpha},$
$X$ carries the \textbf{projective diffeology} $\p$ which is the
pull-back of the diffeologies $\p_{\alpha}$ of each $X_{\alpha}$
via the family of maps $(f_{\alpha})_{\alpha\in\Lambda}.$ The diffeology
$\p$ made of plots $g:O\rightarrow X$ such that, for each $\alpha\in\Lambda,$
\[
f_{\alpha}\circ g\in\p_{\alpha}.\]
 This is the biggest diffeology for which the maps $f_{\alpha}$ are
smooth.

$\bullet$ Let $\lbrace(X_{\alpha},\F_{\alpha},\C_{\alpha})\rbrace_{\alpha\in\Lambda}$
be a family of Fr\"olicher spaces indexed by $\Lambda$ totally ordered
for inclusion,  with $(i_{\beta,\alpha})_{(\alpha,\beta)\in\Lambda^{2}}$
a family of differentiable maps . with natural diffeologies $\p_{\alpha}.$
There is a natural structure of Fr\"olicher space $X=\bigcap_{\alpha\in\Lambda}X_{\alpha},$
which contours \[
\C=\bigcap_{\alpha\in\Lambda}\C_{\alpha}\]
 are the 1-plots of $\p=\bigcap_{\alpha\in\Lambda}\p_{\alpha}.$ A
generating set of functions for this Fr\"olicher space is the set of
maps of the type: 
$$
\bigcup_{\alpha \in \Lambda} \F_{\alpha} \circ f_{\alpha}.$$

\end{Proposition}

\subsection{Differential forms on a diffeological space and differential dimension.}

\begin{Definition} \cite{Sou}
Let $(X,\p)$ be a diffeological space and let $V$ be a vector space equipped with a differentiable structure. 
A $V-$valued $n-$differential form $\alpha$ on $X$ (noted $\alpha \in \Omega^n(X,V))$ is a map 
$$ \alpha : \{p:O_p\rightarrow X\} \in \p \mapsto \alpha_p \in \Omega^n(p;V)$$
such that 

$\bullet$ Let $x\in X.$ $\forall p,p'\in \p$ such that $x\in Im(p)\cap Im(p')$, 
the forms $\alpha_p$ and $\alpha_{p'}$ are of the same order $n.$ 

$\bullet$ Moreover, let $y\in O_p$ and $y'\in O_{p'}.$ If $(X_1,...,X_n)$ are $n$ germs of paths in 
$Im(p)\cap Im(p'),$ if there exists two systems of $n-$vectors $(Y_1,...,Y_n)\in (T_yO_p)^n$ and 
$(Y'_1,...,Y'_n)\in (T_{y'}O_{p'})^n,$ if $p_*(Y_1,...,Y_n)=p'_*(Y'_1,...,Y'_n)=(X_1,...,X_n),$
$$ \alpha_p(Y_1,...,Y_n) = \alpha_{p'}(Y'_1,...,Y'_n).$$

We note by $$\Omega(X;V)=\oplus_{n\in \mathbb{N}} \Omega^n(X,V)$$ the set of $V-$valued differential forms.  
\end{Definition}

With such a definition, we feel the need to make two remarks for the reader:
 
$\bullet$ If there does not exist $n$ linearly independent vectors $(Y_1,...,Y_n)$
defined as in the last point of the definition, $\alpha_p = 0$ at $y.$

$\bullet$ Let $(\alpha, p, p') \in \Omega(X,V)\times \p^2.$ 
If there exists $g \in C^\infty(D(p); D(p'))$ (in the usual sense) 
such that $p' \circ g = p,$ then $\alpha_p = g^*\alpha_{p'}.$ 

\vskip 12pt
\begin{Proposition}
The set $\p(\Omega^n(X,V))$ made of maps $q:x \mapsto \alpha(x)$ from an open subset $O_q$ of a 
finite dimensional vector space to $\Omega^n(X,V)$ such that for each $p \in \p,$
 $$\{ x \mapsto \alpha_p(x) \} \in C^\infty(O_q, \Omega^n(O_p,V)),$$
is a diffeology on $\Omega^n(X,V).$  
\end{Proposition}

Working on plots of the diffeology, one can define the product and the differential of differential forms, which have the same properties as  
the product and the differential of differential forms. 

\begin{Definition}
Let $(X,\p)$ be a diffeological space. 

\noindent
$\bullet$ $(X,\p)$ is \textbf{finite-dimensional} at $x$ if and only if $$\exists n_0\in\mathbb{N},\quad \forall n\in \mathbb{N}, \quad n\geq n_0 \Rightarrow dim(\Omega^n(X,\mathbb{R}))=0$$ 
Then, we set $$dim(X,\p)=max\{n\in \mathbb{N}| dim(\Omega^n(X,\mathbb{R}))>0\}.$$
\noindent
$\bullet$ If not, $(X,\p)$ is called \textbf{infinite dimensional}.
\end{Definition}

Let us make a few remarks on this definition.
If $X$ is a manifold with $dim(X)=n,$ the natural diffeology as described in section \ref{1.1} 
(also called ``n\'ebuleuse'' diffeology) is such that $$dim(X,\p_0)=n.$$
Now, if $(X,\F,\C)$ is the natural Fr\"olicher structure on $X,$ take $\p_1$ generated by the maps of the type
$g\circ c$, where $c\in \C$ and $g$ is a smooth map from an open subset of a finite dimensional space to $\mathbb{R}.$
This is an easy exercise to show that $$dim(X,\p_1)=1.$$
This first point shows that the dimension depends on the diffeology considered.
This leads to the following definition, since $\p(\F)$ is clearly the diffeology with the biggest dimension associated to $(X,\F,\C)$:

\begin{Definition}
The \textbf{dimension} of a Fr\"olicher space $(X,\F,\C)$ is the dimension of the diffeologial space $(X,\p(\F)).$
\end{Definition}

\end{document}